\documentclass[12pt]{article}

\usepackage{amsmath}
\usepackage{amssymb}
\usepackage{latexsym}
\usepackage{amsfonts}
\usepackage{fullpage}
\usepackage{verbatim}

\begin{document}
\centerline{  \hfil Mean Values of Certain Multiplicative Functions \hfil}
\centerline {\hfil and Artin's Conjecture on Primitive Roots  \hfil }

\medskip

\centerline { \hfil Sankar Sitaraman \hfil }

\centerline { \hfil Howard University, Washington, DC \hfil }

\medskip

\centerline {\bf Abstract}

\smallskip

We discuss how one could study  asymptotics of cyclotomic quantities via the mean values of certain multiplicative functions and their Dirichlet series using a theorem of Delange.  We show how this could provide a new approach to Artin's conjecture on primitive roots. We focus on whether a fixed prime has a certain order modulo infinitely many other primes. We also give an estimate for the mean value of one such Dirichlet series.

\bigskip

\centerline {\bf 1. Introduction}

\medskip

Our initial motivation was to come up with a Dirichlet series that would help in obtaining information about the asymptotic behavior of invariants of cyclotomic fields such as class numbers.

 Given a prime $p$ and $\displaystyle \zeta_p = e^{2\pi i /p},$ let $\chi_p : G_p \to \mathbb C^*$ be a character of a group associated with the $p-$th cyclotomic field $K_p = \mathbb Q(\zeta_p).$ For instance $G_p$ could be the Galois group $Gal(K_p/\mathbb Q)$ or one of its subgroups such as the inertia group or the decomposition group. Other possibilities are the ideal class group or its subgroups. 

Consider the Dirichlet series $C(s,\chi) $ defined via an Euler product in the following way:
$$\displaystyle C(s,\chi) = \sum_{m=1}^{\infty} \frac {\chi(m)}{m^s}  = \prod_p \left ( 1 - \frac {\chi_p(g_p)}{p^s} \right )^{-1}$$   
Here $g_p \in G_p$ is chosen suitably.  The multiplicative function $\chi(m)$ is defined by letting $\chi(p) = \chi_p(g_p)$ and then extending to all integers $m$ multiplicatively. Notice that unlike the classical Dirichlet $L-$function $L(s,\chi)$  we use different characters  for different primes $ p.$ This ``dual" or ``horizontal" Dirichlet series seems to be capable of providing information about the asymptotic distribution of certain quantities associated with the $G_p.$ One could probably get much deeper results if $\chi$ were to be considered as a character of a Galois object such as the Galois group of the  compositum of all the $K_p$ (as it is it could be considered a character of $\mathbb Q^*$) but in this paper we simply look at $\chi$ as a multiplicative function.

Given a multiplicative function $f$ let $M_n(f)$ denote the mean value $ \frac 1 n \sum_{m \le n} f(m)$ and $M(f)$ denote the limiting mean value $\displaystyle \lim_{n \to \infty} M_n(f)$ where it exists. 

In Section 3 we provide an estimate for the mean value $M_n(\chi)$ of $\chi$ in $C(s, \chi):$ 

\noindent {\bf Theorem 1}  {\it Given a positive integer $n$ with $\pi(n) =    N,$ where $\pi(n)$ is the number of primes in the set of primes $\{ p_1, p_2,..., p_N\} \subset [1,n],$  $\displaystyle \chi(p_i) = e^{\frac {2\pi i}{d_i} }$ for some positive integer $d_i$ and $i = 1,2,3,....$ chosen such that $d_i \le d_j$ if $i \le j,$  we have}
$$ |M_n(\chi)| = \left | \frac 1 n \sum_{m \le n} \chi(m) \right | \le 
\frac  1 n \sum_{k=1}^{N-1} \frac {d_k}{(N-k)!} \left (  \frac {(log \ n+\sum_{i=k+1}^Nlog \ p_i)^{N-k} }{ \prod_{i=k+1}^N log \ p_i}   \right). \eqno {(1)} $$

 In Section 4 we discuss the applications to asymptotics of cyclotomic quantities that can be obtained by using Delange's theorem on mean values of multiplicative functions. We do not use any specific estimate for mean values, but rather discuss what one can obtain if one had the desired estimates.

The cyclotomic quantity whose asymptotics we discuss are $f_q(p),$ defined below. In the framework for $C(s, \chi)$ above we let $g_p$  be a generator of the group $\displaystyle G_p = Gal(K_p/\mathbb Q)/<\sigma_q>$ where $\sigma_q$ is the Frobenius automorphism of $K_p$ corresponding to the prime $q.$ In the place of $\chi$ we have the multiplicative function $\psi$ given by $\displaystyle \psi(p) = \psi_p(g_p) = e^{2\pi i\frac {f_q(p)}{p-1} } = e^{\frac{2\pi i}{t_q(p)} }$ where $f_q(p)$ is the multiplicative order of the prime $q$ modulo the prime $p.$ We have that $f_q(p)$ is also the order of the Frobenius automorphism $\sigma_q.$  Let  $\psi(\pm 1) = \pm 1.$ The prime $q$ will remain fixed throughout while we consider the asymptotics of $f_q(p)$ over all primes $p.$

Notice also that $\displaystyle t_q(p) = \frac {p-1} {f_q(p)}$ can be defined separately as the number of distinct primes over $q$ in $K_p$ or the number of irreducible factors of the $p-$th cyclotomic polynomial over $\mathbb F_q.$ 

Another example is $\displaystyle C(s, \varpi) $ where $\varpi (p) =  \varpi_p(\sigma_q) = e^{\frac {2\pi i}{f_q(p)} } $
if $p \ne q$ and equals 0 if $p=q.$ In this case $g_p \in  Gal(K_p/\mathbb Q) $ is the Frobenius automorphism $\sigma_q.$ 

We discuss how Dirichlet series such as $C(s,\psi)$ and $C(s,\varpi)$ could be useful in computing asymptotics of $f_q(p).$ 
We prove the following, mainly using Delange's theorem:

\noindent {\bf Theorem 2} {\it   Let $\varpi(p) = e^{\frac {2\pi i }{f_q(p)} }$ if $p \ne q$ and equals 0 if $p=q$ where $f_q(p)$ is as defined above.  

\smallskip

(a) The mean value of $\varpi \to c$ for some nonzero $c$ and $C(s,\varpi)$ has residue $c$ at $s=1.$

\smallskip

(b) If $\psi$ is as defined above and the mean value of $\psi $ denoted by $M(\psi)$ equals zero then there are infinitely many $p$ such that $\displaystyle f_q(p) > \frac {p-1}{log \ p}.$}

\medskip

  Much work has been done on the asymptotics of $f_q(p).$  Artin conjectured that $ f_a(p) = p-1$ infinitely often with non-zero density $C$ (Artin constant), for any non-zero integer $a$ other than 1,-1 or a perfect square. Progress has been made on Artin's conjecture by Hooley \cite{Hooley},Gupta-Murty \cite{Gupta-Murty}, Heath-Brown \cite{Heath-Brown} and others. Heath-Brown proved that Artin's conjecture is true for one of 2, 3 and 5.

 In Section 4 also  we provide a procedure for proving that, given any prime $q,$ there are infinitely many primes $p$ for which the order $f_q(p) $ takes values in a certain set. [Although we restrict $q$ to be a prime, this procedure should work for other integers as well]. We also provide a formula for the Artin density for any prime $p$ using  another Dirichlet series. 

  We believe that this is a new approach to the problem of asymptotics of cyclotomic quantities. Washington \cite{Washington} laid the groundwork and Friedman \cite{Friedman} studied invariants of the product of Iwasawa extensions of cyclotomic fields over several primes. In a different direction,  Davenport-Heilbronn \cite{Davenport-Heilbronn} and Shintani \cite{Shintani} constructed objects that compute the asymptotic distribution of class numbers of cubic fields. The Cohen-Lenstra heuristics  also throw some light on the asymptotics of class numbers. The work on Artin's conjecture for primitive roots can also be considered as related to asymptotics of $f_q(p)$ but much of the previous work has been about $f_q(p)$ for a family of primes $q.$ Our focus here is on whether a fixed prime $q$ has a certain order modulo infinitely many other primes. 

I thank Francois Ramaroson and Niranjan Ramachandran for useful conversations and Chris McCarthy for help with computations in the initial stages of this work. I also thank the anonymous referee for his valuable comments.

\medskip  

\centerline {\bf 2.  Preliminaries }

\medskip

We use the following theorem by Delange to relate the mean value of a multiplicative function to  the values at primes. This is but one of several such theorems proved by Erdos-Wintner, Halasz and others.  A nice exposition is given in Tenenbaum \cite{Tenenbaum}. Perhaps one could even use the Selberg-Delange method to get deeper results.

\noindent {\bf Theorem} [Delange] {\it  Let $g$ be a multiplicative function with values in the unit disc.  If the limiting mean value $\displaystyle M(g) := \lim_{ x \to \infty} \frac 1 x \sum_{n \le x} g(n)$  is non-zero then $\displaystyle \sum_p \frac {1-g(p)}p < \infty.$ Moreover if  $\displaystyle \sum_p \frac {1-g(p)}p < \infty$ then $g$ has a limiting mean value that is given by  $\displaystyle \prod_{p} (1-p^{-1}) \sum_{\nu = 0}^{\infty}g(p^\nu)p^{-\nu}.  $}

Also, we see that if $g$ is completely (or strictly) multiplicative, then $\displaystyle \sum_p \frac {1-g(p)}p < \infty$ iff  $\displaystyle \prod_{p } (1-p^{-1}) \sum_{\nu = 0}^{\infty}g(p^\nu)p^{-\nu}  =  \prod_{p} \frac{1-p^{-1}}{1-g(p)p^{-1}}  $ is non-zero. 
So in this case by Delange's theorem the existence of a non-zero limiting mean value is equivalent to the convergence of $\displaystyle \sum_p \frac {1-g(p)}p.$

\medskip

\centerline  {\bf  3. Estimate for mean value}

\medskip
 As before,  $p$ is a prime and $\chi_p : G_p \to \mathbb C^*$ is a character of a group associated with the $p-$th cyclotomic field $K_p = \mathbb Q(\zeta_p).$ In most of the cases below we will denote by $g_p$ a generator of the corresponding cyclic group $G_p.$

It is easy to see that the following Dirichlet series $C(s,\chi) $ defined by its Euler product is convergent for
$Re(s) > 1$ and that the Euler product expansion is valid:
$$\displaystyle C(s,\chi) = \sum_{m=1}^{\infty} \frac {\chi(m)}{m^s}  = \prod_p \left ( 1 - \frac {\chi_p(g_p)}{p^s} \right )^{-1}$$   

  The multiplicative function $\chi(m)$ is defined by letting $\chi(p) = \chi_p(g_p)$ and extending to all $m$ multiplicatively.

An estimate for the mean value of $C(s, \chi)$ is useful for determining when it has analytic continuation to $s < 1$ and to get information about the behavior at $s=1.$ To get an idea about how the mean value of $C(s,\chi)$ behaves we looked at $\displaystyle  C_0(s) = \prod_p \left ( 1 - \frac {\zeta_p}{p^s} \right )^{-1}.$ This series may not be useful in computing asymptotics of cyclotomic quantities and it was merely used as the motivation for Theorem 1 below. It does not explicitly appear in the discussion below. Nevertheless we believe that $C_0(s)$ would be useful in problems such as the computation of the number of smooth integers and the question of uniform distribution of certain arithmetic functions.

Instead of relying on existing methods for character sums or exponential sums we have used a simpler approach that we consider better suited for the properties peculiar to $\chi(m).$ 

For a given ordering of the primes  $p_i$  we write $\displaystyle \chi(p_i) = e^{\frac {2\pi i}{d_i} }$ for some positive integer $d_i$ and $i = 1,2,3,....$

Let $\displaystyle m = \prod_i p_i^{k_i}$ be the prime decomposition of the positive integer $m.$

We write the mean value as 
$$ M_n(\chi) = \frac 1 n \sum_{m \le n} \chi(m)
= \frac 1 n \sum_{m \le n} e^{ 2 \pi i \left ( \sum \frac {k_i}{d_i}  \right )} 
= \frac 1 n \sum_{m \le n} \prod_{p_i | m} e^{ 2 \pi i \frac {k _i}{d_i} } .$$

To get the estimate we use the cancellations among the various terms that results from the basic fact that
$1+ \zeta_d + \zeta_d^2 + ....+ \zeta_d^{d-1} = 0$ for any positive integer $d$ and $\displaystyle \zeta_d = e^{\frac{2 \pi i}{d} }.$ 
We need to look at
the space of the exponents $k_i$ and reduce them modulo the $d_i.$ Suppose that 
$d_1 \le d_2 \le d_3 \le ....\le d_N$ then let the corresponding primes less than $n$ be $p_1,p_2,....,p_N.$ Note that {\it the ordering is based on that of the $d_i$ and not of the primes.}  If $m$ corresponds to 
$(k_1,k_2,...,k_N) \in \mathbb Z^N$ (some of the exponents could be zero) then the terms corresponding to $(k_1,k_2,...,k_i+1,...,k_N), (k_1,k_2,...,k_i+2,...,k_N),....,(k_1,k_2,...,k_i+d_i,...,k_N)$
together add up to zero. This is the idea we use to prove the following:

\noindent {\bf Proof of Theorem 1.} \quad We will give a geometric proof based on the idea described above. 

We start by looking at 
the contributions from $\chi(m)$ for all $m \le n$ such that $p_1 | m,$ i.e, $\displaystyle \sum_{p_1 | m} \chi(m).$  The set $S_1 = \{ ( k_1,...,k_N) \  | \ k_1 \ne 0, \ \sum_{i=1}^N k_i log\ p_i \le log \ n \}$ is the ``exponent space" of the $m$  in this sum.  In this set $S_1$ we can remove all the elements in a given ``stalk" by which we mean all the elements ``above" a given $(0,k_2,....,k_N)$ which have points in ``strings" of length $d_1,$ namely sets of points of the form 
$(rd_1+1,k_2,...,k_i,...,k_N), (rd_1+2, ,k_2,...,k_i,...,k_N),....,((r+1)d_1, k_1,k_2,...,k_i,...,k_N)$ for $r$ a nonnegative integer. Thus ``reducing" this ``stalk" modulo $d_1$ we get a smaller set or the ``tip" of the stalk with a set of points of length less than $d_1.$ Doing this for each $(0,k_2,...,k_N)$ we get above each point of the ``hyperpolyhedron" $H_1 = \{ (k_1,k_2,....,k_N) \   | \  k_1 =0, \sum_{i=2}^N k_ilog\ p_i \le log \ n \}$ a truncated ``stalk" of length less than or equal to $d_1.$ The number of lattice points in this ``reduced" $S_1$ which we call $\bar {S_1}$ is less than the volume of the polyhedral ``box" with base $ \{ (0,k_2,...,k_N) \in H_1\  | \ \sum_{i=2}^N ( k_i-1) log\ p_i \le log \ n \}$ and height $d_1.$ 

Notice that in the last line we are not using the polyhedron bounded by the condition $\sum k_i log \ p_i \le log \ n $ that comes from $m \le n.$ We are using a bigger polyhedron whose volume will be bigger than  the volume of the $N-$dimensional object obtained by enclosing each lattice point by a ``box" with length 1 in the positive direction in every dimension. The latter is a standard argument for counting lattice points. [cf., for example, Granville \cite{Granville}].

Hence we get  (with $m \in \bar{S_1}$ meaning $(k_1,k_2,...,k_N) \in \bar{S_1} $ by abuse of notation) 

$$ \left | \sum_{p_1 | m} \chi(m) \right | \le \sum_{m \in \bar{S_1} } |\chi(m) | \le  
\frac {d_1}{(N-1)!} \left ( \frac {(log \ n+\sum_{i=2}^Nlog \ p_i)^{N-1} }{ \prod_{i=2}^N log \ p_i}\right) .$$

Next we consider the contribution from $m$ such that $(p_1,m) = 1, p_2|m$ or equivalently, corresponding to points with $k_1=0, k_2 \ne 0.$ Notice that in the previous step we did not remove any $\chi(m)$ corresponding to $m$ with exponent co-ordinates of the form $(0,k_2,...,k_N).$ The ``exponent space" of such $m$ which we shall call $S_2$ can be ``reduced" in a similar fashion until  the truncated space $\bar {S_2}$ is arrived at.  All the steps carried out for $S_1$ will go through for the space $S_2$ and the contribution from the terms with $m \in S_2$ will be 

$$ \left | \sum_{p_1\nmid  m, p_2 |m} \chi(m) \right | \le \sum_{m \in \bar{S_2} } |\chi(m) | \le  
\frac {d_2}{(N-2)!} \left (\frac {(log \ n+\sum_{i=3}^Nlog \ p_i)^{N-2} }{ \prod_{i=3}^N log \ p_i}  \right) .$$

Proceeding inductively and adding up the contributions from $m $ in $S_1,S_2,...., S_{N-1}$ we get the inequality (1). {\it QED}

\medskip

\noindent {\it Remark: } The above estimate could possibly be sharpened by using estimates for character sums after reducing the sum modulo the relation $1+ \zeta_d + \zeta_d^2 + ....+ \zeta_d^{d-1} = 0.$ One could also replace the ``boxes" with smaller polyhedra that also enclose the lattice points. 

\bigskip

\centerline   {\bf 4. $\mathbf{C(s,\psi)}$ and Artin's conjecture on primitive roots}

\bigskip

We discuss how the mean value of Dirichlet series of the form $C(s,\chi)$ could be useful in studying the asymptotics of cyclotomic quantities. Specifically we look at the cyclotomic quantity $f_q(p)$ defined, as in the Introduction, as the multiplicative order of the prime $q$ modulo the prime $p$. Also $\displaystyle t_q(p) = \frac {p-1} {f_q(p)}.$ The two Dirichlet series $C(s, \psi)$ and $\displaystyle C(s, \varpi) $ are also as defined in the Introduction. Their Euler products run over all primes $p$ while the prime $q$ remains fixed.

We look at two cases  namely $|M_n(\chi)| \to 0 $ and $|M_n(\chi) | \to c$ (for $\chi = \psi \ or \ \varpi$) where $0 < c \le 1$ is a fixed real number. As is well known $|M_n(\chi) | \to c$  implies that the residue of $C(s, \chi)$ at $s=1$  is $c.$  In this paper we do not seek to prove that $|M_n(\chi)| \to 0 $ or $|M_n(\chi) | \to c$ for these Dirichlet series but rather look at the consequences. Nevertheless we indicate how this could possibly be done,  in the following two brief remarks:

\medskip

\noindent {\it Remarks: }

1. It is possible that one could prove these using the estimate of  Theorem \ref{mean-estimate} above in combination with other estimates, especially if the $d_i$ (in the case of $\psi, d_i = t_q(p_i))$ remain small. The crucial ingredient one needs is an estimate for $\chi(m)$ in terms of $m.$ Good estimates can probably be obtained if we frame the $\chi(m)$ (such as $\psi(m)$ or $\varpi(m)$) in the context of a Galois object such as the Galois group of the  compositum of all the $K_p.$ Another possibility is to relate $m$ to the $f_q(p)$ for $p | m$ using a formula such as 
 $\displaystyle i_q(m) = \sum_{d | m} \frac {\phi(d)}{f_q(d)}.$  Here $i_q(m)$ is the number of distinct irreducible factors of $x^m-1$ over $\mathbb F_q.$

2. On the other hand one could try to estimate $|M(\chi)|$ by separating the primes. More precisely, if it is known that the sum of $\chi(m)$ are bounded in a nice way for $m$ divisible by  certain primes which contribute the most, then using summation by parts one could write $M_n(\chi)$ in terms of such $\displaystyle  \chi(m)$ and then get the estimate for $|M(\chi)|,$ assuming there are ways to easily estimate the contribution from the other primes.

\medskip

\medskip

\noindent {\bf Proof of Theorem 2. } 

\noindent (a) By Delange's theorem the mean value has a nonzero limit if
$\displaystyle \sum_p \frac {1-\varpi(p)} p < \infty.$   Now trivially for any $p \ne q$ we have $f_q(p) > k \ 
log \ p.$ for a suitable constant $k.$  From this it is easy to show that 
$\displaystyle  Re(1-e^{\frac {2\pi i }{f_q(p)} } )< c_1/(log\ p)^2$ and 
$\displaystyle -Im(1-e^{\frac {2\pi i }{f_q(p)} } ) < c_2/log\ p$ for some constants $c_1,c_2$ and all but a finite number of primes. 
From this it follows that $\displaystyle \sum_p \frac {1-\varpi (p)} p < \infty$ since both $\displaystyle \sum_p \frac 1 {p(log\ p)}$ and $\displaystyle \sum_p \frac 1 {p(log\ p)^2}$ are finite. Hence the mean value of $\varpi \to c$ for some nonzero $c$ and $C(s,\varpi)$ has residue $c$ at $s=1.$

\smallskip

\noindent {\it Remark: }
 In the case when $q=2$ by
 Erdos-Murty \cite{erdos-murty} it is known that for almost all primes $p$ we have $\displaystyle f_2(p) > p^{1/2 + \epsilon(p))}$ for any
 $\epsilon(p) \to 0.$

\medskip

\noindent (b) Also by Delange's theorem if the mean value goes to zero then
$\displaystyle \sum_p \frac {1-\psi(p)} p = \infty.$   If it is true that only finitely many $p$ have $\displaystyle f_q(p) > \frac {p-1}{log\ p} $  then for the rest of the primes we get 
$\displaystyle f_q(p) \le \frac {p-1} {log\ p} < \frac p {log \ p}. $ This  implies $\displaystyle \frac{f_q(p)}{p-1} < \frac 2 {log\ p} $ for all but finitely many primes. From this we get that $\displaystyle \sum_p \frac {1-\psi (p)} p < \infty$ in the same way as in (a) which contradicts the hypothesis. \quad {\it QED}

\smallskip

\noindent {\it Remark: }
Numerical evidence and probabilistic arguments support the conclusion that $M(\varpi) $ is nonzero 
and $M(\psi) = 0.$ [Look at the sums over primes in the RHS in Delange's theorem].

\medskip

\bigskip

Next we outline a procedure to relate $C(s,\psi)$ to the values of $f_q(p).$ We first demonstrate it with the Dirichlet $L-$function $L(s, \chi)$ where $\chi$ is a Dirichlet character modulo a prime $q.$ Note that $\chi$ is completely multiplicative. Since the character group is cyclic in this case, we can take $\chi$ to be a generator. Then the density of primes with a given value for $\chi(p)$ is found using the orthogonality of $\chi^k$ for $k = 1,2,...,q-1$ and the fact that $L(1,\chi^k) \ne 0$ for $k = 1,2,...,q-2.$

In the case of the $C(s,\psi)$ defined above we don't know of an orthogonality relation, and also it is much harder to prove results of the type $L(1,\chi^k) \ne 0.$ So we needed to see what can be obtained without them. To understand this, let us see what can be done with $L(1,\chi^k) $ without orthogonality or $L(1,\chi^k) \ne 0.$

We have $\displaystyle  | \sum_{n \le x} \chi^k(n) |  \le q-1$ for $k = 1,2,...,q-2.$ So in this case the limiting mean value $M(\chi^k) = 0$ for all $k = 1,2,...,q-2.$ By Delange's theorem (and the remarks following it) we get that $\displaystyle \sum_p \frac {1-\chi^k(p)}p $ does not converge, for $k = 1,2,...,q-2.$ Now $1-\chi^k(p) = 0$ whenever $k$ is divisible by the order of $\chi(p)$ in $\mathbb C^*.$
Let $S_k = \{ p \ prime \ | \ o(\chi(p)) \nmid \ k  \}$ where $o(\chi(p)) $ is the order of $\chi(p)$ in $\mathbb C^*.$ Then 
$\displaystyle \sum_{p \in S_k} \frac {1-\chi^k(p)}p $ does not converge for $k = 1,2,...,q-2$ means that
 $S_k$ must be infinite for $k = 1,2,...,q-2.$ 

\medskip
  Let $M(\psi^k)$ be the limiting mean value of $\psi^k.$ We get the following theorem, whose proof relies on an argument very similar to the one outlined above for $\chi.$

\noindent {\bf Theorem 3} {\it  Let $\displaystyle t_p = \frac {p-1}{f_q(p)}.$ If $M(\psi^k)$ does not exist or is zero, then $S_k = \{ p \ prime \ | \ t_p \nmid \ k  \}$ is infinite.}

\bigskip

\noindent {\it \bf Artin Density:} 

\medskip

Now we give a formula to compute the Artin density for $q,$ namely $\displaystyle \lim_{x \to \infty} N_q(x)/\pi(x)$ where $N_q(x)$ is the number primes $p \le x$  such that $q$ is a primitive root mod $p.$ 

In fact we give a formula for 
$  \displaystyle  \sum_{ t_p = t} p^{-s}$ from which the corresponding Dirichlet density for $t$ can be found. (If $t=1$ we get the Artin density for $q$). If the natural density exists then the Dirichlet density also exists, so we hope that this is of interest.

We use the following multiplicative function $\psi_t:$ 

$\psi_t(p) = -1$ if 
$\displaystyle \frac {p-1}{f_q(p)} = t$ and $\psi_t(p) = 1$ otherwise, where 
$t \in \mathbb Z^{>0}.$ 

As we did for $\psi$ we let $\psi_t(\pm 1) = \pm 1.$ Then extend $\psi_t$   to all of $\mathbb Z$ multiplicatively. 

With $\zeta(s)$ being the Riemann zeta function, expand the Euler products of $\zeta(s)$ and $C(s, \psi_t)$ in terms of the powers of the primes and then look at $ log\  \zeta (s) -  log\ C(s, \psi_t) .$ The contribution from the terms involving $p^{ks}$ for $k \ge 2$ is bounded. In the remaining part $\displaystyle \sum_p \frac {1- \psi_t(p) } {p^s}$ those terms containing $ 1- \psi_t(p)$ for $p$ such that  $\displaystyle \frac {p-1}{f_q(p)} \ne t$  vanish and we get the following, for $s > 1$ (where everything on the RHS converges):

$$ 2 \sum_{ t_p = t} p^{-s} = log\  \zeta (s) -  log\ C(s, \psi_t) + O(1)  \eqno{(2)}$$

The reason we think this formula is  worth mentioning is the following. Note that the limit as $s \to 1$ of the right hand side is basically the log of the limiting mean value $M(\psi_t),$ if it exists, according to Delange's theorem. And there maybe some way to estimate this mean value, since it depends on values at $n,$ although we will somehow have to get around the dependence on knowing which $n$ are divisible by $p$ such that $t_p = t.$

\bigskip

\centerline {\bf 5. Further questions }

\medskip

\smallskip

1. What are the analytic properties of $C(s,\psi),$ namely the regions of convergence and analytic continuation, zeroes, functional equation, special values, etc ? The problems discussed above could all be framed in terms of the special values of these Dirichlet series.

2. What are the algebraic properties of $C(s,\psi),$ namely its relation to the structure of the Galois group of the composita of the $K_p$ and computation of  cyclotomic invariants?

3. Is there any relationship to modular forms? For instance, the Euler product of $C(s,\psi)C(s, \overline{\psi})$ looks like that of a modular form of weight 1 with coefficients$\ a_p$ given by $a_p =  \psi(p) + \overline{\psi(p)}.$ 

4. While we did not have it in mind while working on this, we wonder if the ``pretentious" approach to analytic number theory developed by Granville and Soundararajan (cf. for instance, Balog-Granville-Soundararajan \cite{pretentious}) might be useful here, in place of or in addition to Delange's theorem?

\end{document}